\documentclass{elsart3-1}

\usepackage{graphicx,amsmath,amssymb}  
\def\ZZ{\mathbb Z}

\def\bx{{\bf x}}
\def\by{{\bf y}}
\def\bz{{\bf z}}

\def\bk{{\bf k}}

\def\br{{\bf r}}

\def\btau{\boldsymbol{\tau}}

\def\b1{{\bf 1}}
\def\d1{\mathds{ 1}}

\usepackage[english,francais]{babel}

\newtheorem{e-proposition}[theorem]{Proposition}

\newtheorem{e-definition}[theorem]{Definition\rm}

\renewcommand{\theequation}{\arabic{equation}}
\def\og{\leavevmode\raise.3ex\hbox{$\scriptscriptstyle\langle\!\langle$~}}
\def\fg{\leavevmode\raise.3ex\hbox{~$\!\scriptscriptstyle\,\rangle\!\rangle$}}

\journal{the Acad\'emie des sciences}
\begin{document}
\centerline{}
\begin{frontmatter}
%



%
\selectlanguage{english}
\title{On the lower bound of the discrepancy  of Halton's sequence}



\author[authorlabel1]{Mordechay B. Levin},
\ead{mlevin@math.biu.ac.il}

\address[authorlabel1]{Department of Mathematics,
Bar-Ilan University, Ramat-Gan, 52900, Israel}



\medskip
\selectlanguage{francais}
\begin{center}
{\small Re\c{c}u le *****~; accept\'e apr\`es r\'evision le +++++\\
Pr\'esent\'e par £££££}
\end{center}

\selectlanguage{english}
\begin{abstract}
\selectlanguage{english}
\vskip 0.5\baselineskip
Let $ (H_s(n))_{n \geq 1} $ be an $s-$dimensional Halton's sequence. Let $D_N$ be the discrepancy of the sequence $ (H_s(n))_{n = 1}^{N} $.
It is known that $ND_N =O(\ln^s N)$ as $N \to  \infty $. 
 In this paper we prove that this estimate  is exact: 
 $$
      \overline{\lim}_{ N \to \infty}   N \ln^{-s}(N) D_N >0.
 $$
{\it To cite this article: A. Nom1,  C. R.
Acad. Sci. Paris, Ser. I .}

\selectlanguage{francais}
\noindent{\bf R\'esum\'e} \vskip 0.5\baselineskip \noindent
{\bf Sur la limite inf\'erieure de la discr\'epance de la s\'equence de Halton. }\\
Soit $ (H_s (n))_{n \geq 1} $  une s\'equence de Halton $ s- $ dimensions. Soit $ D_N $  la discr\'epance de la s\'equence $ (H_s (n))_{n = 1}^{N} $. 
Il est connu que $N D_N = O (\ln^s N) $ que $ N \to \infty $. 
  Dans cet article, nous montrons que cette estimation  est exacte: 
 $$
      \overline{\lim}_{ N \to \infty}   N \ln^{-s}(N) D_N >0.
 $$
{\it Pour citer cet article~: A. Nom1,  C. R.
Acad. Sci. Paris, Ser. I .}
\vskip 0.5\baselineskip
\end{abstract}
\end{frontmatter}

\selectlanguage{english}



\textheight=20.8 cm
\section{Introduction}
 Let $(\beta_{n})_{n \geq 1}$ be a  sequence in  the unit cube $[0,1)^s$, $\by =(y_1,...,y_s)$,
 $B(\by)=[0,y_1) \times \cdots \times [0,y_s) \subseteq [0,1)^s $,
\begin{equation}\label{1}
\Delta(B(\by), (\beta_{n})_{n=1}^{N}  )= \sum_{1 \leq n \leq N}  ( \b1_{B(\by)}(\beta_{n}) -  y_1 \cdots y_s), \quad {\rm where} \quad  \b1_{B(\by)}(\bx) =1, \; {\rm if} \; 
\bx  \in B(\by), 
\end{equation}
and $   \b1_{B(\by)}(\bx) =0,$  if $ 
\bx   \notin B(\by)$. We define the star {\it discrepancy} of a 
$N$-point set $(\beta_{n})_{n=1}^{N}$ as
\begin{equation} \label{2}
   \emph{D}^{*}((\beta_{n})_{n=1}^{N}) = 
    \sup_{ 0<y_1, \ldots , y_s \leq 1} \; | \frac{1}{  N}
  \Delta(B(\by),(\beta_{n})_{n=1}^{N}) |.
\end{equation}

Let $(\beta_{n})_{n \geq 1}$ be an arbitrary sequence in  $[0,1)^s$.
In 1954, Roth proved that 
\begin{equation} \nonumber
   \limsup_{N \to \infty } N (\ln N)^{-\frac{s}{ 2}} \emph{D}^{*}((\beta_{n})_{n=1}^{N})>0 . 
\end{equation}
According to the well known conjecture (see, e.g., [1, p.283]), this estimate can be improved
\begin{equation}   \label{4}
 \limsup_{N \to \infty } N (\ln N)^{-s} \emph{D}^{*}((\beta_{n})_{n=1}^{N})>0 . 
\end{equation}
In 1972, W. Schmidt proved this conjecture for $ s=1 $. For $s=2$, Faure and Chaix [4] proved (\ref{4}) for a class of $(t,s)-$sequences.
For a review of research on this conjecture, see for example  [2].\\

{ \bf Definition.} {\it  An $s$-dimensional sequence  $(\beta_{n})_{n \geq 1}$ is of 
 \texttt{low discrepancy} (abbreviated
l.d.s.) if \\$ \emph{D}^{*}((\beta_{n})_{n=1}^{N})=O(N^{-1}(\ln
N)^{s}) $ for $ N \rightarrow \infty $. }

Let $ p\geq 2 $ be an integer,
 \begin{equation}\label{1.4}
 n=\sum_{i\geq 0}e_{p,i}(n) p^i,\; {\rm with} \; e_{p,i}(n) \in \{0,1, \ldots
 ,p-1\}, \qquad  {\rm and}    \qquad  \phi_p(n)= \sum_{i\geq 0}e_{p,i}(n) p^{-i-1}.
 \end{equation}
Van der Corput (see [3, ref. 1891])   proved that $ (\phi_p(n))_{n\geq 0}$ is a $1-$dimensional l.d.s. Let
\begin{equation}\label{1.5} 
  H_s(n)= (\phi_{p_1}(n),\ldots ,\phi_{p_s}(n)), \quad n=0,1,2,...,
\end{equation}
 where $
p_1,\ldots ,p_s\geq 2 $ are pairwise coprime integers.
Halton (see [3, ref. 729])   proved that $ ( H_s(n))_{n\geq 0}$ is an $s-$dimensional l.d.s. 
 For other examples of  l.d.s. see e.g.  [1], [3].
  In $\mathsection 2$  we will prove \\

{\bf Theorem.} {\it Let $p_0=p_1p_2\cdots p_s$, $s \geq 2$ and $m_0=[2 p_0 \log_2 p_0] +2$. Then }
\begin{equation}\nonumber
      \sup_{1 \leq N \leq 2^{mm_0}}   N \emph{D}^{*}(( H_s(n))_{n=1}^{N}) \geq m^s(8p_0)^{-1} \quad {\rm for} \quad m \geq p_0. 
\end{equation}  
{\bf Remark.} This result  supports the conjecture (\ref{4}).  In [5], we received  a similar result for  sequences obtained from algebraic lattices. 
In [6], we proved  a similar result for  some $(t,s)-$sequences (see also [7]).

\section{Proof of the Theorem.}
Let $x_i =\sum_{j \geq 1}  x_{i,j}p_i^{-j}$, with $x_{i,j} \in \{0,1,...,p_i-1  \}$, $[x_i]_r =\sum_{1 \leq j \leq r}  x_{i,j}p_i^{-j}$, $i=1,...,s$, $r=1,2,...$\\
By (\ref{1.4}), we have  $\phi_{p_i}(n) \in [[x_i]_r, [x_i]_r +p_i^{-r})$ if and only if $n \equiv   \dot{x}_{i,r}  \;\; ({\rm mod} \; p_i^r)$, where $\dot{x}_{i,r} = \sum_{1 \leq j \leq r}  x_{i,j}p_i^{j-1}$.
Let $\br =(r_1,...,r_s)$,  $P_{\br} = p_1^{r_1} \cdots p_s^{r_s}$ and   $M_{i,\br} \equiv (P_{\br}  p_i^{-r_i})^{-1}   \; ({\rm mod} \; p_i^{r_i})$.
Using  the Chinese Remainder Theorem, we get
\begin{equation}\label{1.12}
    \phi_{p_i}(n) \in [[x_i]_{r_i}, [x_i]_{r_i} +p_i^{-r_i}),\;\; {\rm for} \;\; i=1,...,s  \;   \Longleftrightarrow   n \equiv   \ddot{x}_{\br}  \; ({\rm mod} \; P_{\br})
 \; {\rm with} 	\;	\ddot{x}_{\br} =\sum_{i=1}^s  M_{i,\br}   P_{\br} p_i^{-r_i}  \dot{x}_{i,r_i}.
\end{equation}
It is easy to verify that if $r_i^{'} \geq r_i$, for all $i=1,...,s$, then
\begin{equation} \label{1.14}
              \ddot{x}_{\br^{'} } \equiv \ddot{x}_{\br}  \;\; ({\rm mod} \; P_{\br}).
\end{equation}
We consider the case $x_{i,r_i} \neq 0 $ for all $i=1,...,s$. We obtain from (\ref{1.12}) that 
\begin{equation}\label{1.12a}
   \phi_{p_i}(n) \in [[x_i]_{r_i} -p_i^{-r_i}, [x_i]_{r_i} ),\quad {\rm for} \quad i=1,...,s  \;   \Longleftrightarrow   n \equiv   \ddot{x}_{\br} - \sum_{i=1}^s  M_{i,\br}   P_{\br} p_i^{-1} \;\; ({\rm mod} \; P_{\br}).
\end{equation}
Let $p_0 =p_1 p_2\cdots p_s$, $\breve{p}_i =p_0/p_i$, $\tau_i =\min \{ 1 \leq k <\breve{p}_i | p_i^k \equiv 1 \;\; ({\rm mod} \; \breve{p}_i)  \}$, $i=1,...,s$.
Let $\by =(y_1,...,y_s)$ with $y_i = \sum_{1 \leq j \leq m} p_i^{-j\tau_i}$,   $[y_i]_{\tau_ik_i}=\sum_{1 \leq j \leq k_i} p_i^{-j\tau_i} $,  and let
$\dot{y}_{i,\tau_ik_i}=\sum_{1 \leq j \leq k_i} p_i^{j\tau_i -1} $, 
$k_i \geq 1$,  $i=1,...,s$,  $ \bk=(k_1,...,k_s)$,
   $B(\by)=[0,y_1) \times \cdots \times [0,y_s) \subset [0,1)^s $,  $B_{\bk} =\prod_{1 \leq i \leq s} [[y_i]_{\tau_ik_i} -p_i^{-k_i\tau_i},[y_i]_{\tau_ik_i}  )$,\\ $\btau = (\tau_1,...,\tau_s)$,  $u \cdot v =(u_1v_1,...,u_sv_s)$.
We have
\begin{equation}\label{1.16}
    B(\by) =\bigcup_{1 \leq k_1,...,k_s \leq m} B_{\bk}, 		\quad {\rm and} \quad 
		\b1_{B(\by)} (\bz) -  y_1 \cdots y_s = 		
		\sum_{1 \leq k_1,...,k_s \leq m}    ( \b1_{B_{\bk}} (\bz) -  P_{\btau \cdot \bk}^{-1}).
\end{equation} 
Let  $\hat{y}_{\btau \cdot \bk} =\sum_{i=1}^s  M_{i,\btau \cdot \bk}, 
 P_{\btau \cdot \bk} p_i^{-\tau_i k_i} \dot{y}_{i,\tau_i(k_i-1)}$
 and 
\begin{equation}\label{1.19}
   A_{\bk} \equiv -\sum_{i=1}^s  M_{i,\btau \cdot \bk}   P_{\btau \cdot \bk} p_i^{-1}  \;\; ({\rm mod} \; P_{\btau \cdot \bk}),   \quad {\rm with} \quad  A_{\bk} \in[0,  P_{\btau \cdot \bk}) . 
\end{equation}
From (\ref{1.12})  we get  $\ddot{y}_{\btau \cdot \bk} = \sum_{i=1}^s  M_{i,\btau \cdot \bk}, 
 P_{\btau \cdot \bk} p_i^{-\tau_i k_i} \dot{y}_{i,\tau_ik_i} \equiv \hat{y}_{\btau \cdot \bk}  - A_{\bk} \;\; ({\rm mod} \; P_{\btau \cdot \bk} )$.  
 By (\ref{1.5}) and (\ref{1.12a}), we obtain
\begin{equation}\label{1.18}
     H_s(n) \in B_{\bk}  \;   \Longleftrightarrow    \phi_{p_i}(n) \in [\dot{y}_{i,\tau_i k_i} -p_i^{-\tau_i k_i}, \dot{y}_{i,\tau_i k_i} ),\quad {\rm for} \quad i=1,...,s     \Longleftrightarrow   n \equiv   \hat{y}_{\btau \cdot \bk}  \;\; ({\rm mod} \; P_{\btau \cdot \bk}). 
\end{equation}
Let
\begin{equation}\label{1.19a}
 \tilde{y}_m :\equiv \hat{y}_{\btau(m+1) } \; ({\rm mod} \; P_{\btau(m+1) }),  \; {\rm with} \; \tilde{y}_{ m }  \in[0,  P_{\btau(m+1)}), \; {\rm where} \; 
  \btau(m+1) =(\tau_1 (m+1),...,\tau_s (m+1)).
\end{equation}
Using (\ref{1.14}),  we get 
 $    \hat{y}_{\btau \cdot \bk} - A_{\bk}   \equiv \ddot{y}_{\btau \cdot \bk}  \equiv
  \ddot{y}_{\btau(m+1)} \equiv \hat{y}_{\btau(m+1)}  - A_{\btau(m+1)}  \equiv  \tilde{y}_m \;\; ({\rm mod} \; P_{\btau \cdot \bk})$, with $k_1,...,k_s$\\$ \in [1,m]$.
Applying (\ref{1.18}), we have
\begin{equation}  \nonumber
     H_s(n) \in B_{\bk}  \;   \Longleftrightarrow   n \equiv   \tilde{y}_m +A_{\bk}
      \;\; ({\rm mod} \; P_{\btau \cdot \bk}).  
\end{equation}
 By (\ref{1.18}), we get
\begin{equation} \label{1.22}  
        \sum_{n =  \tilde{y}_m + N_1 P_{\btau \cdot \bk} }^{ \tilde{y}_m +(N_1+1)P_{\btau \cdot \bk}-1}  ( \b1_{B_{\bk}} (H_s(n)) - P_{\btau \cdot \bk}^{-1} )  =0,
				\quad  {\rm and} \qquad 
				 \sum_{n= \tilde{y}_m + N_1 P_{\btau \cdot \bk}}^{ \tilde{y}_m +N_1P_{\btau \cdot \bk}+N_2-1}  ( \b1_{B_{\bk}} (H_s(n)) - P_{\btau \cdot \bk}^{-1} ) 
\end{equation}
\begin{equation} \nonumber
 =\sum_{ n \in [\tilde{y}_m, \tilde{y}_m +N_2)}  ( \b1_{B_{\bk}} (H_s(n)) - P_{\btau \cdot \bk}^{-1} )
				= \sum_{ \substack{ n \in [\tilde{y}_m, \tilde{y}_m +N_2) \\ n= \tilde{y}_m +A_{\bk} }}  1 - N_2P_{\btau \cdot \bk}^{-1}  
       =  \b1_{[0,N_2)}(A_{\bk}) -N_2P_{\btau \cdot \bk}^{-1},
\end{equation}
with $N_1 \geq 0$ and $N_2 \in [0, P_{\btau \cdot \bk})$, $N_1,N_2 \in \ZZ$. 
From (\ref{1}) and (\ref{1.16}), we get
\begin{equation}\nonumber
        \Delta(B(\by), ( H_s(n))_{n=\tilde{y}_m}^{\tilde{y}_m+N-1} ) 
  =       \sum_{y_{0,m} \leq n < y_{0,m} +N}  ( \b1_{B(\by)}(H_s(n)) -  y_1 \cdots y_s)  
\end{equation}
\begin{equation}\label{1.32}
      = \sum_{1 \leq k_1,...,k_s \leq m} \rho(\bk,N), \quad {\rm with} \quad  
 \rho(\bk,N)= 
       \sum_{y_{0,m} \leq n < y_{0,m} +N}  ( \b1_{B_{\bk}} (H_s(n)) -  P_{\btau \cdot \bk}^{-1})  .
\end{equation}
Let
\begin{equation} \label{1.32a}
\alpha_m:=   \frac{1}{P_{\btau m}} \sum_{ N=1}^{ P_{\btau m}} 
  \Delta(B(\by), ( H_s(n))_{n=\tilde{y}_m}^{\tilde{y}_m+N-1} )   = \sum_{1 \leq k_1,...,k_s \leq m} \alpha_{m,\bk}, \;\; {\rm with} \;\;  
 \alpha_{m,\bk}= \frac{1}{P_{\btau m}} \sum_{ N=1}^{ P_{\btau m}} \rho(\bk,N).
\end{equation}
Bearing in mind  (\ref{1.22}) and (\ref{1.32}) , we derive
\begin{equation}\nonumber
  \alpha_{m,\bk} =
      \frac{1}{P_{\btau m}} \sum_{ N_1=0}^{ P_{\btau m}/ P_{\btau \cdot \bk} -1}\sum_{ N_2=1}^{ P_{\btau \cdot \bk}}  \Big(
       \sum_{ n=\tilde{y}_m}^{\tilde{y}_m  +N_1P_{\btau \cdot \bk} -1}  ( \b1_{B_{\bk}} (H_s(n)) -  P_{\btau \cdot \bk}^{-1})  
\end{equation}
\begin{equation}\nonumber
    +   \sum_{ n=\tilde{y}_m  +N_1P_{\btau \cdot \bk}}^{\tilde{y}_m  +N_1P_{\btau \cdot \bk}+N_2-1}   ( \b1_{B_{\bk}} (H_s(n)) -  P_{\btau \cdot \bk}^{-1})  \Big)
     =  
      \frac{1}{P_{\btau m}} \sum_{ N_1=0}^{ P_{\btau m}/ P_{\btau \cdot \bk} -1}\sum_{ N_2=1}^{ P_{\btau \cdot \bk}}  \Big(    \b1_{[0,N_2)}(A_{\bk}) -N_2P_{\btau \cdot \bk}^{-1}   \Big) 
\end{equation}
\begin{equation}\nonumber
        = \frac{1}{P_{\btau \cdot \bk}} \sum_{ N_2=1}^{ P_{\btau \cdot \bk}}  \Big(    \b1_{[0,N_2)}(A_{\bk}) -N_2P_{\btau \cdot \bk}^{-1}   \Big) = \frac{P_{\btau \cdot \bk} -A_{\bk}}{P_{\btau \cdot \bk}}-\frac{P_{\btau \cdot \bk}(P_{\btau \cdot \bk}+1)}{2P_{\btau \cdot \bk}^2} = \frac{1}{2} -   \frac{A_{\bk}}{P_{\btau \cdot \bk}}  - \frac{1}{2P_{\btau \cdot \bk}}.
\end{equation}
Using  (\ref{1.32a}), we have
\begin{equation}  \label{1.23a}
\alpha_m  =  \sum_{1 \leq k_1,...,k_s \leq m}  \Big(\frac{1}{2} -   \frac{A_{\bk}}{P_{\btau \cdot \bk}}  - \frac{1}{2P_{\btau \cdot \bk}} \Big).
\end{equation}
Taking into account that  $M_{i,\btau \cdot \bk} \equiv (P_{\btau \cdot \bk}  p_i^{-\tau_i k_i})^{-1}  \equiv \prod_{1 \leq j \leq s, j \neq i} 
   p_j^{-\tau_j k_j}  \;\; ({\rm mod} \; p_i^{\tau_i k_i})$, and that \\
$p_j^{\tau_j} \equiv 1 \;\; ({\rm mod} \; p_i) \; (i \neq j) $,     we obtain   $ M_{i,\btau \cdot \bk}  \equiv 1 \;\; ({\rm mod} \; p_i)$, $i=1,...,s$. 
From (\ref{1.19}), we get 
\begin{equation}\nonumber
  [0,1) \ni \frac{A_{\bk}}{P_{\btau \cdot \bk}} \equiv  -  \sum_{1 \leq i \leq s}  M_{i,\btau \cdot \bk}   P_{\btau \cdot \bk} p_i^{-1}/P_{\btau \cdot \bk}  \equiv -\frac{1}{p_1} - \cdots  - \frac{1}{p_s}  \;\; ({\rm mod} \; 1).
\end{equation}
Applying  (\ref{1.23a}), we derive
\begin{equation}\label{1.23}
\alpha_m  
 = m^s \Big( \frac{1}{2} 
      - \{-\beta  \} \Big) - 
      \sum_{1 \leq k_1,...,k_s \leq m} \frac{1}{2P_{\btau \cdot \bk}},  \quad {\rm with} \quad  \beta =  \frac{1}{p_1} + \cdots  + \frac{1}{p_s},
\end{equation}
where $\{x\}$ is the fractional part of  $x$.
Let $ \beta    \equiv 1/2 \;\; ({\rm mod} \; 1)$. Hence $p_0 =p_1p_2\cdots p_s \equiv 0 \;\; ({\rm mod} \; 2)$. 
Let $p_{\nu} \equiv 0 \;\; ({\rm mod} \; 2)$ for some $\nu \in [1,s]$. Then 
\begin{equation} \nonumber
    b_1:=    p_0(p_{\nu}/2-1)/p_{\nu}              \equiv p_0(\beta -1/p_{\nu}) = p_0\sum_{1 \leq i \leq s, i \neq \nu}  1/p_i  ({\rm mod} \; 1)  \; 
          {\rm and} \;  b_1 \equiv b_2   ({\rm mod} \; p_0), \; {\rm with}\;
					b_2 = \sum_{i \neq \nu}  p_0/p_i.	
\end{equation}
Let $j \in [1,s]$ and $j \neq \nu$. We see that $b_1  \equiv 0 \;\; ({\rm mod} \; p_j)$ and $b_2  \not \equiv 0 \;\; ({\rm mod} \; p_j)$.
We get a contradiction. Hence  $ \beta   \not \equiv 1/2 \;\; ({\rm mod} \; 1)$. We have
\begin{equation}\nonumber
 0 \neq \Big| \frac{1}{2} - \Big\{  - \beta\Big\} \Big| =
     \Big| \frac{1}{2} - \Big\{-   \Big(\frac{1}{p_1} + \cdots  + \frac{1}{p_s} \Big) \Big\} \Big| =  \frac{|a|}{2p_0}, \quad {\rm with \;some \;integer} \; a. 
\end{equation}
Thus $  | 1/2 - \{  - \beta \} |  \geq 1/(2p_0)   $. 
Bearing in mind that $P_{\btau \cdot \bk} \geq 2^{k_1+k_2 + \cdots +k_s}  $, we obtain from (\ref{1.23})
\begin{equation}\label{1.30}
|\alpha_m| \geq   
 \frac{m^s}{2p_0}  -\frac{1}{2} =  \frac{m^s}{2p_0} (1 - \frac{p_0}{ m^{s}}) \geq \frac{m^s}{4p_0}
     \quad {\rm for} \quad m \geq p_0>4.
\end{equation}
It is easy to see that $\tau_i \leq p_0$, $i=1,...,s$, and $ 2P_{\btau (m+1)} = 2 p_1^{\tau_1 (m+1)} \cdots p_s^{\tau_s (m+1)} \leq 2^{1+p_0 (m+1) \log_2 p_0} $ \\$\leq 2^{m(1+
2p_0  \log_2 p_0)}$
$\leq 2^{mm_0}$ with  $m_0=[2 p_0 \log_2 p_0] +2$. 
Using  (\ref{1.19a}), we have that $\tilde{y}_m + P_{\btau m} <2P_{\btau (m+1)} $ \\$  \leq 2^{mm_0}$.
By  (\ref{1.30}), (\ref{1.32a}) and (\ref{2}),  we get
\begin{equation}\nonumber
  m^s (4p_0)^{-1} \leq  |\alpha_m |   \leq  \sup_{1  \leq N \leq P_{\btau m}}   N \emph{D}^{*}(( H_s(n))_{n=\tilde{y}_m}^{\tilde{y}_m+N-1})     
\end{equation}  
\begin{equation}\nonumber
  \leq  \sup_{1 \leq L,L+N \leq 2P_{\btau (m+1)}}   N \emph{D}^{*}(( H_s(n))_{n=L}^{L+N-1})   \leq  2\sup_{1 \leq N \leq 2^{\dot{mm_0}}}   N \emph{D}^{*}(( H_s(n))_{n=1}^{N}) \quad {\rm for} \quad m\geq p_0.
\end{equation} 
Hence, the Theorem is proved.

\end{document}